\documentclass[12pt,oneside]{amsproc}
\usepackage[koi8-r]{inputenc}
\usepackage[russian]{babel}
\usepackage{amssymb,epic,eepic}
\numberwithin{equation}{section}
\newtheorem{lem}{Лемма}[section]
\newtheorem{tm}{Теорема}[section]

\title[]{Индефинитная задача Штурма--Лиувилля для некоторых классов
самоподобных сингулярных весов}
\author{А.~А.~Владимиров,%
\address{Московский государственный университет
им.~М.~В.~Ломоносова, механико-математический факультет}
\email{vladimi@mech.math.msu.su}
И.~А.~Шейпак%
\address{Московский государственный университет
им.~М.~В.~Ломоносова, механико-математический факультет}
\email{iasheip@mech.math.msu.su}}
\thanks{Работа поддержана РФФИ, грант \No\,04-01-00712, и фондом поддержки
ведущих научных школ, грант~НШ-1927.2003.1.}
\newcommand{\Wo}{{\raisebox{0.2ex}{\(\stackrel{\circ}{W}\)}}{}}
\newcommand{\ind}{\operatorname{ind}}

\begin{document}
\noindent УДК~517.984
\begin{abstract}
В статье продолжается изучение вопроса об асимптотике спектра граничной задачи
\begin{gather*}
	-y''-\lambda\rho y=0,\\ y(0)=y(1)=0,
\end{gather*}
где \(\rho\) есть функция из пространства \(\Wo_2^{-1}[0,1]\), имеющая
самоподобную первообразную. Рассматриваются случаи неарифметического и
вырожденного арифметического самоподобия такой первообразной.
\end{abstract}
\maketitle

\section{Введение}
В настоящей статье будет продолжено начатое в работе~\cite{VlSh} изучение
вопроса о спектральных асимптотиках граничной задачи
\begin{gather}\label{eq:0.1}
	-y''-\lambda\rho y=0,\\ \label{eq:0.2} y(0)=y(1)=0,
\end{gather}
где весовая функция \(\rho\) имеет самоподобную квадратично суммируемую
первообразную. Цель настоящей статьи состоит в рассмотрении двух случаев
самоподобия, оставшихся неразобранными в предыдущей работе.

Первым из них является возникающий при допущении индефинитности веса \(\rho\)
вырожденный случай арифметического самоподобия, описываемый ниже в
теореме~\ref{tm0}. Соответствующий аналог теоремы восстановления будет доказан в
параграфе~\ref{pt:1}.

Вторым из этих случаев является случай неарифметического самоподобия.
Соответствующий аналог теоремы восстановления будет доказан в
параграфе~\ref{pt:2}.

\section{Теорема восстановления в вырожденном арифметическом
случае}\label{pt:1}
Основным результатом настоящего параграфа является следующее утверждение.
\begin{tm}\label{tm:1.1}
Пусть фиксировано произвольное положительное вещественное число \(\tau\).
Пусть наборы неотрицательных чисел \(\{u_{k}\}_{k=1}^{N}\) и \(\{v_k\}_{k=1}^N\)
таковы, что
\[
	\sum\limits_{k=1}^N (u_{k}+v_{k})=1,\qquad
	\sum\limits_{k=1}^N v_k>0,
\]
причём наибольший общий делитель номеров \(k\), для которых справедливо
неравенство \(u_k+v_k>0\), равен \(1\). Пусть также
\begin{align*}
	&\forall k\leqslant [(N-1)/2]&&u_{2k+1}=0,\\
	&\forall k\leqslant [N/2]&&v_{2k}=0,
\end{align*}
где через \([a]\) обозначается целая часть числа \(a\). Наконец, пусть
непрерывные на \(\mathbb R\) функции \(X_j\), где \(j=1,2\), удовлетворяют
условиям
\begin{align*}
	&\forall t\in\mathbb R^+\qquad |X_j(t)|\leqslant \Pi_j\,e^{-\tau t},\\
	&\forall t\in\mathbb R^-\qquad X_j(t)=0.
\end{align*}
Здесь \(\Pi_j\), где \(j=1,2\), "--- положительные числа. Тогда существует и
единственна пара непрерывных на \(\mathbb R\) функций \(Z_j\), где \(j=1,2\),
удовлетворяющая уравнениям
\begin{align*}
	&\forall t\in\mathbb R^+\qquad Z_j(t)=X_j(t)+\sum\limits_{k=1}^N (u_k
	\,Z_j(t-k)+v_k\,Z_{3-j}(t-k)),\\
	&\forall t\in\mathbb R^-\qquad Z_j(t)=0.
\end{align*}
При этом для функций \(Z_j\), где \(j=1,2\), справедливы оценки
\[
	\forall t\in\mathbb R^+\qquad \left|Z_j(t)-s(t-j+1)\right|\leqslant
	(\Pi_1+\Pi_2)\cdot C(t),
\]
где исчезающая на бесконечности функция \(C\) не зависит от выбора
функций \(X_j\), а непрерывная \mbox{\(2\)-пе}\-ри\-о\-ди\-чес\-кая функция
\(s\) имеет вид
\[
	\forall t\in\mathbb R\qquad s(t)=\dfrac{1}{J}
	\sum\limits_{k=-\infty}^{+\infty}\left(X_1(t-2k)+X_2(t-2k-1)\right).
\]
Через \(J\) здесь обозначена величина
\[
	J\rightleftharpoons\sum\limits_{k=1}^N k\,(u_k+v_k).
\]

\end{tm}

Доказательство теоремы~\ref{tm:1.1} полностью аналогично доказательству
теоремы~2.2 из работы~\cite{VlSh}. Поэтому мы не станем приводить его подробного
изложения, а укажем лишь на лежащую в его основе лемму. При этом, как и
в~\cite{VlSh}, через \(\ell_{1,r}\), где \(r\in\mathbb R^+\), мы будем
обозначать пространство последовательностей \(\theta=\{\theta_n\}_{%
n=0}^{\infty}\), ограниченных по норме
\[
	\|\theta\|_{\ell_{1,r}}\rightleftharpoons
	\sum\limits_{n=0}^{\infty}r^n|\theta_n|.
\]
\begin{lem}\label{lemA:1}
Пусть фиксированы произвольные \(r<1\) и \(R>1\). Пусть наборы неотрицательных
чисел \(\{u_{k}\}_{k=1}^{N}\) и \(\{v_{k}\}_{k=1}^{N}\) таковы, что
\[
	\sum\limits_{k=1}^N (u_{k}+v_{k})=1,\qquad
	\sum\limits_{k=1}^N v_{k}>0,
\]
причём наибольший общий делитель номеров \(k\), для которых справедливо
неравенство \(u_k+v_k>0\), равен \(1\). Пусть также
\begin{align*}
	&\forall k\leqslant [(N-1)/2]&&u_{2k+1}=0,\\
	&\forall k\leqslant [N/2]&&v_{2k}=0,
\end{align*}
где через \([a]\) обозначается целая часть числа \(a\). Тогда для любых
\(x_j\in\ell_{1,R}\), где \(j=1,2\), существует и единственна пара
\(z_j\in\ell_{1,r}\) решений системы
\[
	\forall n\in\{0,1,2,\ldots\}\qquad z_{j,n}=x_{j,n}+
	\sum\limits_{k=1}^{\min(N,n)}(u_kz_{j,n-k}+v_kz_{3-j,n-k}),
	\qquad j=1,2.
\]
При этом координаты векторов \(z_j\), где \(j=1,2\), удовлетворяют оценке
\[
	\forall n\in\{0,1,2,\ldots\}\qquad\left|z_{j,n}-\dfrac{\omega-
	(-1)^{n+j}\chi}{J}\right|\leqslant(\|x_1\|_{\ell_{1,R}}+
	\|x_2\|_{\ell_{1,R}})\cdot C_{n}, \qquad j=1,2,
\]
где бесконечно малая последовательность \(\{C_k\}_{k=0}^{\infty}\) не
зависит от выбора \(x_j\), а величины \(\omega\), \(\chi\) и \(J\)
определены равенствами
\begin{gather*}
	\omega\rightleftharpoons\dfrac12
	\sum\limits_{k=0}^{\infty}(x_{1,k}+x_{2,k}),\\
	\chi\rightleftharpoons\dfrac12
	\sum\limits_{k=0}^{\infty}(-1)^k(x_{1,k}-x_{2,k}),\\
	J\rightleftharpoons\sum\limits_{k=1}^N k(u_k+v_k).
\end{gather*}
\end{lem}
\begin{proof}
Для доказательства леммы достаточно повторить, с незначительными изменениями,
рассуждения из доказательства леммы~2.2 работы~\cite{VlSh}. Поэтому мы
ограничимся указанием на основное отличие рассматриваемого случая от
разобранного в~\cite{VlSh}, опуская технические детали. При этом, как и
в~\cite{VlSh}, через \(X_j\), \(Z_j\), \(U\) и \(V\) мы будем обозначать
производящие функции последовательностей \(\{x_{j,n}\}_{n=0}^{\infty}\),
\(\{z_{j,n}\}_{n=0}^{\infty}\), \(\{u_n\}_{n=1}^{N}\) и \(\{v_n\}_{n=1}^{N}\).

В рассматриваемом вырожденном случае производящие функции \(Z_j\), где
\(j=1,2\), имеют на границе единичного круга не одну особенность (в точке
\(w=1\)), а две (в точках \(w=\pm 1\)). Учёт дополнительной особенности легко
проводится, если принять во внимание, что при выполнении условий леммы
справедливо тождество
\[
	\forall w\in\mathbb C\qquad (1-U+V)(w)=(1-U-V)(-w).
\]
Оно позволяет разложить величины \((1-U+V)(w)\) на множители \(1+w\) и
\(Q(-w)\), где \(Q\) "--- многочлен, удовлетворяющий тождеству
\[
	\forall w\in\mathbb C\qquad (1-U-V)(w)=(1-w)\cdot Q(w).
\]
Сказанное означает, что последовательности \(\{z_{j,n}\}_{n=0}^{\infty}\)
асимптотически эквивалентны последовательностям коэффициентов Маклорена функций
\[
	\dfrac{X_1(1)+X_2(1)}{2J(1-w)}+\dfrac{X_j(-1)-X_{3-j}(-1)}{2J(1+w)}.
\]
Отсюда и вытекает утверждение доказываемой леммы.
\end{proof}

\section{Теоремы восстановления в неарифметическом случае}\label{pt:2}
Целью настоящего параграфа является доказательство следующих утверждений.
\begin{tm}\label{tmA:1}
Пусть набор положительных чисел \(\{u_{k}\}_{k=1}^{N}\) таков, что
\[
	\sum\limits_{k=1}^N u_{k}=1.
\]
Пусть также набор положительных чисел \(\{l_k\}_{k=1}^N\) линейно независим над
полем \(\mathbb Q\). Наконец, пусть непрерывная на \(\mathbb R\) функция \(X\)
удовлетворяет при \(t\to\pm\infty\) условию
\[
	X(t)=o(t^{-2}).
\]
Тогда существует и единственна непрерывная на \(\mathbb R\) функция \(Z\),
удовлетворяющая уравнениям
\begin{gather*}
	\forall t\in\mathbb R\qquad Z(t)=X(t)+\sum\limits_{k=1}^N u_k
	\,Z(t-l_k),\\
	\lim\limits_{t\to-\infty} Z(t)=0.
\end{gather*}
При этом функция \(Z\) подчиняется соотношению
\[
	\lim\limits_{t\to+\infty}Z(t)=\dfrac{1}{J}\int\limits_{-\infty}^{+\infty}
	X\,d\mu,
\]
где через \(d\mu\) обозначена линейная мера Лебега, а через \(J\) обозначена
величина
\[
	J\rightleftharpoons\sum\limits_{k=1}^N u_k\,l_k.
\]
\end{tm}
\begin{tm}\label{tmA:2}
Пусть наборы неотрицательных чисел \(\{u_{k}\}_{k=1}^{N}\) и \(\{v_k\}_{k=1}^N\)
таковы, что
\begin{gather}\notag
	\forall k\in\{1,\ldots,N\}\qquad u_k+v_k>0,\\ \notag
	\sum\limits_{k=1}^N v_k>0,\\ \label{eq3}
	\sum\limits_{k=1}^N (u_{k}+v_{k})=1.
\end{gather}
Пусть также набор положительных чисел \(\{l_k\}_{k=1}^N\) линейно независим над
полем \(\mathbb Q\). Наконец, пусть непрерывные на \(\mathbb R\) функции
\(X_j\), где \(j=1,2\), удовлетворяют при \(t\to\pm\infty\) условиям
\begin{equation}\label{eq:as}
	X_j(t)=o(t^{-2}).
\end{equation}
Тогда существует и единственна пара непрерывных на \(\mathbb R\) функций
\(Z_j\), где \(j=1,2\), удовлетворяющая уравнениям
\begin{gather}\label{eq1}
	\forall t\in\mathbb R\qquad Z_j(t)=X_j(t)+\sum\limits_{k=1}^N (u_k
	\,Z_j(t-l_k)+v_k\,Z_{3-j}(t-l_k)),\\ \label{eq2}
	\lim\limits_{t\to-\infty} Z_j(t)=0.
\end{gather}
При этом функции \(Z_j\) подчиняются соотношениям
\begin{equation}\label{eq5}
	\lim\limits_{t\to+\infty}Z_j(t)=\dfrac{1}{2J}
	\int\limits_{-\infty}^{+\infty} (X_1+X_2)\,d\mu,
\end{equation}
где через \(d\mu\) обозначена линейная мера Лебега, а через \(J\) обозначена
величина
\[
	J\rightleftharpoons\sum\limits_{k=1}^N (u_k+v_k)\,l_k.
\]
\end{tm}

Непосредственно нами будет проведено доказательство теоремы~\ref{tmA:2}.
Теорема~\ref{tmA:1} может быть доказана аналогичным образом; кроме того, она,
по-существу, не нова. В частности, для случая экспоненциального убывания функции
\(X\) утверждения теоремы~\ref{tmA:1} установлены в работе~\cite{LV} (см.
также~\cite[гл.~IX, \S~9]{Fell}).

Доказательство теоремы~\ref{tmA:2} будет опираться на ряд вспомогательных
утверждений.

\begin{lem}\label{lem1}
Пусть непрерывные на \(\mathbb R\) функции \(X_j\), где \(j=1,2\),
неотрицательны, и пусть непрерывные на \(\mathbb R\) функции \(Z_j\), где
\(j=1,2\), удовлетворяют уравнениям~\eqref{eq1},~\eqref{eq2}. Тогда функции
\(Z_j\) также неотрицательны.
\end{lem}
\begin{proof}
Предположим, что при некоторых \(t\in\mathbb R\), \(m\in\{1,2\}\) и
\(\varepsilon>0\) справедливо неравенство \(Z_m(t)<-\varepsilon\). Если бы
при всех \(k\in\{1,\ldots,N\}\) и \(j\in\{1,2\}\) были справедливы неравенства
\(Z_j(t-l_k)\geqslant-\varepsilon\), то, в силу равенства~\eqref{eq3} и
неотрицательности функций \(X_j\), из уравнения~\eqref{eq1} следовало бы
противоречащее сделанному предположению неравенство
\(Z_m(t)\geqslant-\varepsilon\). Следовательно, при некоторых \(k\in
\{1,\ldots,N\}\) и \(j\in\{1,2\}\) должно выполняться неравенство
\(Z_j(t-l_k)<-\varepsilon\). Положительность величин \(l_k\) приводит теперь
к противоречию с уравнением~\eqref{eq2}.
\end{proof}

\begin{lem}\label{lem2}
Существует число \(C>0\) со следующим свойством: для любых
непрерывных на \(\mathbb R\) функций \(X_j\) и \(Z_j\), где \(j=1,2\),
удовлетворяющих уравнениям~\eqref{eq1},~\eqref{eq2} и оценкам
\begin{equation}\label{eq4}
	\forall t\in\mathbb R\qquad |X_j(t)|\leqslant\dfrac{\Pi}{t^2+1},
\end{equation}
где \(\Pi\) "--- произвольное неотрицательное число, справедливы также оценки
\[
	\forall t\in\mathbb R\qquad |Z_j(t)|\leqslant C\cdot\Pi.
\]
\end{lem}
\begin{proof}
\textit{Шаг~1.} Рассмотрим функцию \(\eta\) вида
\[
	\forall t\in\mathbb R\qquad \eta(t)=\sum\limits_{k=1}^N (u_k+v_k)\cdot
	(\operatorname{arctg}(t)-\operatorname{arctg}(t-l_k)).
\]
Эта функция положительна и имеет при \(t\to\pm\infty\) асимптотику
\(\eta\asymp t^{-2}\). Если теперь положить \(\tilde X_1=\tilde X_2
\rightleftharpoons\eta\) и
\[
	\forall t\in\mathbb R\qquad \tilde Z_1(t)=\tilde Z_2(t)\rightleftharpoons
	\dfrac{\pi}{2}+\operatorname{arctg}(t),
\]
то подстановка функций \(\tilde X_j\) и \(\tilde Z_j\) в
уравнения~\eqref{eq1},~\eqref{eq2} вместо, соответственно, функций \(X_j\) и
\(Z_j\) будет давать верные равенства.

\textit{Шаг~2.} Зафиксируем постоянную \(C>0\) со свойством
\[
	\forall t\in\mathbb R\qquad C\geqslant\dfrac{\pi}{\eta(t)\cdot(t^2+1)}.
\]
Тогда для функций \(X_j\), удовлетворяющих условию~\eqref{eq4}, будут
выполняться неравенства
\[
	\forall t\in\mathbb R\qquad |X_j(t)|\leqslant \dfrac{C\cdot\Pi}{\pi}\,
	\eta(t).
\]
Из леммы~\ref{lem1} и результатов предыдущего шага теперь вытекают неравенства
\[
	\forall t\in\mathbb R\qquad |Z_j(t)|\leqslant \dfrac{C\cdot\Pi}{\pi}\,
	\left(\dfrac{\pi}{2}+\operatorname{arctg}(t)\right),
\]
означающие справедливость утверждения леммы.
\end{proof}

\begin{lem}\label{lem5}
Пусть непрерывные на \(\mathbb R\) функции \(X_j\), где \(j=1,2\),
удовлетворяют при \(t\to\pm\infty\) условиям~\eqref{eq:as}. Тогда существует
и единственна пара непрерывных на \(\mathbb R\) функций \(Z_j\), где \(j=1,2\),
удовлетворяющая уравнениям~\eqref{eq1},~\eqref{eq2}.
\end{lem}
\begin{proof}
Единственность искомой пары решений легко выводится из леммы~\ref{lem1}.
Докажем существование этой пары.

Если функции \(X_j\), где \(j=1,2\), тождественно равны нулю левее некоторой
точки \(t_0\), то соответствующие решения \(Z_j\) находятся тривиальным образом:
левее \(t_0\) они полагаются тождественно равными нулю, а правее \(t_0\) их
значения определяются из уравнений~\eqref{eq1}. Заметим теперь, что
если функции \(X_j\) подчиняются асимптотике~\eqref{eq:as}, то найдутся
такие последовательности исчезающих в окрестности точки \(-\infty\) функций
\(\{X_{j,n}\}_{n=1}^{\infty}\), что будут справедливы соотношения
\begin{equation}\label{eq:appr}
	\lim\limits_{n\to\infty}\left(\sup\limits_{t\in\mathbb R}
	\biggl(|X_{j,n}(t)-X_j(t)|\cdot(t^2+1)\biggr)\right)=0.
\end{equation}
Из леммы~\ref{lem2} теперь следует, что последовательности
\(\{Z_{j,n}\}_{n=1}^{\infty}\) решений задач
\begin{gather}\label{eq1:+}
	\forall t\in\mathbb R\qquad Z_{j,n}(t)=X_{j,n}(t)+\sum\limits_{k=1}^N
	(u_k\,Z_{j,n}(t-l_k)+v_k\,Z_{3-j,n}(t-l_k)),\\ \label{eq2:+}
	\lim\limits_{t\to-\infty} Z_{j,n}(t)=0
\end{gather}
равномерно на \(\mathbb R\) сходятся к некоторым решениям \(Z_j\)
задачи~\eqref{eq1},~\eqref{eq2}. Тем самым, лемма доказана.
\end{proof}

\begin{lem}\label{lem4}
Для любых вещественных чисел \(\zeta_j\), где \(j=1,2\), удовлетворяющих
неравенству \(\zeta_1+\zeta_2\neq 0\), найдутся такие финитные непрерывные
функции \(X_j\), где \(j=1,2\), что будут справедливы равенства
\[
	\int\limits_{-\infty}^{+\infty}X_j\,d\mu=\zeta_j,
\]
а соответствующие решения \(Z_j\) задачи~\eqref{eq1},~\eqref{eq2} будут
подчиняться асимптотике~\eqref{eq5}.
\end{lem}
\begin{proof}
Зафиксируем произвольную непрерывную функцию \(Z_1\), тождественно равную \(0\)
в некоторой окрестности \(-\infty\), и тождественно равную \((\zeta_1+\zeta_2)/
2J\) в некоторой окрестности \(+\infty\). Зафиксируем также функцию \(Z_2\)
вида
\[
	\forall t\in\mathbb R\qquad Z_2(t)=Z_1(t-\tau),
\]
где \(\tau\) "--- произвольная вещественная постоянная.

Сопоставленные функциям \(Z_j\) на основе уравнений~\eqref{eq1} функции
\(X_j\), очевидно, являются финитными. При этом, как проверяется
непосредственно, имеют место равенства
\begin{gather*}
	\int\limits_{-\infty}^{+\infty}X_1\,d\mu=\dfrac{\zeta_1+\zeta_2}{2}+
	\dfrac{\tau\,(\zeta_1+\zeta_2)}{2J}\sum\limits_{k=1}^N v_k,\\
	\int\limits_{-\infty}^{+\infty}X_2\,d\mu=\dfrac{\zeta_1+\zeta_2}{2}-
	\dfrac{\tau\,(\zeta_1+\zeta_2)}{2J}\sum\limits_{k=1}^N v_k.
\end{gather*}
Таким образом, если выбрать
\[
	\tau\rightleftharpoons\dfrac{J(\zeta_1-\zeta_2)}{(\zeta_1+\zeta_2)\,
	\sum\limits_{k=1}^N v_k},
\]
то будут выполнены все утверждения леммы.
\end{proof}

\begin{lem}\label{lem6}
Пусть финитные функции \(X_j\in C_0^2(\mathbb R)\), где \(j=1,2\),
удовлетворяют условию
\begin{equation}\label{eq:pas2}
	\int\limits_{-\infty}^{+\infty}(X_1+X_2)\,d\mu=0.
\end{equation}
Тогда решения \(Z_j\) уравнений~\eqref{eq1},~\eqref{eq2} подчиняются асимптотике
\[
	\lim\limits_{t\to+\infty}Z_{j}(t)=0.
\]
\end{lem}
\begin{proof}
\textit{Шаг~1.}
Поскольку функции \(X_j\) финитны, их фурье-образы \(\hat X_j\) вида
\[
	\forall t\in\mathbb R\qquad \hat X_j(t)=\dfrac{1}{\sqrt{2\pi}}
	\int\limits_{-\infty}^{+\infty}e^{-it\cdot} X_j\,d\mu
\]
являются бесконечно дифференцируемыми. При этом, ввиду гладкости функций
\(X_j\), фурье-образы \(\hat X_j\) и все их производные имеют в окрестности
точек \(\pm\infty\) асимптотику \(O(t^{-2})\). Кроме того, из
предположения~\eqref{eq:pas2} вытекает равенство \(\hat X_1(0)+\hat X_2(0)=0\).

Зафиксируем теперь произвольную функцию \(\varphi\in C_0^3(\mathbb R)\),
тождественно равную \(1\) в некоторой окрестности точки \(0\), и рассмотрим
последовательности функций \(\{\hat X_{j,n}\}_{n=1}^{\infty}\) вида
\[
	\forall t\in\mathbb R\qquad \hat X_{j,n}(t)=\varphi(t/n)\cdot\hat X_j(t).
\]
Как несложно показать, последовательности \(\{\hat X_{j,n}\}_{n=1}^{\infty}\)
и \(\{\hat X''_{j,n}\}_{n=1}^{\infty}\) сходятся в пространстве
\(L_1(\mathbb R)\) к функциям \(\hat X_j\) и \(\hat X''_j\), соответственно.
Тем самым, для фурье-прообразов \(X_{j,n}\) функций \(\hat X_{j,n}\) будут
выполняться соотношения~\eqref{eq:appr}.

\textit{Шаг~2.} Обозначим через \(U\) и \(V\) квазимногочлены вида
\begin{gather*}
	\forall t\in\mathbb R\qquad U(t)=\sum\limits_{k=1}^N u_k e^{-il_k\,t},\\
	\forall t\in\mathbb R\qquad V(t)=\sum\limits_{k=1}^N v_k e^{-il_k\,t}.
\end{gather*}
Ввиду линейной независимости величин \(\{l_k\}_{k=1}^N\) над полем
\(\mathbb Q\), квазимногочлен \(1-U-V\) имеет на \(\mathbb R\) единственный
(причём простой) нуль в точке \(0\), а квазимногочлен \(1-U+V\) не имеет на
\(\mathbb R\) ни одного нуля.

Рассмотрим теперь последовательности \(\{\hat Z_{j,n}\}_{n=1}^{\infty}\)
финитных функций класса \(C_0^2(\mathbb R)\), удовлетворяющих системам уравнений
\begin{gather*}
	(1-U-V)(\hat Z_{1,n}+\hat Z_{2,n})=\hat X_{1,n}+\hat X_{2,n},\\
	(1-U+V)(\hat Z_{1,n}-\hat Z_{2,n})=\hat X_{1,n}-\hat X_{2,n}.
\end{gather*}
Существование таких функций гарантировано равенствами \(\hat X_{1,n}(0)+
\hat X_{2,n}(0)=0\) и характером гладкости функции \(\varphi\). Несложно
заметить, что фурье-прообразы \(Z_{j,n}\) функций \(\hat Z_{j,n}\) будут
удовлетворять уравнениям~\eqref{eq1:+},~\eqref{eq2:+}. Из результатов
предыдущего шага и леммы~\ref{lem2} (применяемой отдельно к вещественным
и мнимым частям функций \(Z_{j,n}\) и \(X_{j,n}\)) вытекает теперь, что
последовательности \(\{Z_{j,n}\}_{n=1}^{\infty}\) равномерно на \(\mathbb R\)
сходятся к решениям уравнений~\eqref{eq1},~\eqref{eq2}. Однако при всех
\(j=1,2\) и \(n=1,2,\ldots\) заведомо справедливы равенства
\[
	\lim\limits_{t\to+\infty}Z_{j,n}(t)=0.
\]
Тем самым, лемма доказана.
\end{proof}

\begin{proof}[Доказательство теоремы~\ref{tmA:2}]
Зафиксируем произвольное число \(\varepsilon>0\). Зафиксируем также числа
\(\zeta_j\), где \(j=1,2\), удовлетворяющие неравенствам
\(\zeta_1+\zeta_2\neq 0\) и
\[
	\left|\int\limits_{-\infty}^{+\infty}X_j\,d\mu-\zeta_j\right|<
	\varepsilon.
\]
Поскольку функции \(X_j\), где \(j=1,2\), удовлетворяют
асимптотике~\eqref{eq:as}, они могут быть представлены в виде
\[
	X_j=X_j^0+X_j^1+X_j^2,
\]
где:
\begin{itemize}
\item функции \(X_j^0\in C_0^2(\mathbb R)\) удовлетворяют равенству
\[
	\int\limits_{-\infty}^{+\infty}(X_1^0+X_2^0)\,d\mu=0;
\]
\item функции \(X_j^1\) сопоставлены величинам \(\zeta_j\) на основе
леммы~\ref{lem4};
\item функции \(X_j^2\) подчиняются неравенствам
\[
	\forall t\in\mathbb R\qquad |X_j^2(t)|\leqslant
	\dfrac{\varepsilon}{t^2+1}.
\]
\end{itemize}
Комбинируя утверждения лемм~\ref{lem2}, \ref{lem4} и~\ref{lem6}, устанавливаем,
что для функций \(Z_j\) асимптотически при \(t\to+\infty\) справедливы
оценки
\[
	\left|Z_j-\dfrac{1}{2J}\int\limits_{-\infty}^{+\infty}(X_1+X_2)\,d\mu
	\right|\leqslant\varepsilon\cdot\left(C+\dfrac{1}{J}\right),
\]
где постоянная \(C\) определена леммой~\ref{lem2}. Поскольку \(\varepsilon\)
выбрано произвольным образом, то тем самым теорема доказана.
\end{proof}

\section{Асимптотика спектра}
Перейдём теперь к рассмотрению спектральных свойств операторного пучка
\(T_{\rho}:\Wo_2^1[0,1]\to\Wo_2^{-1}[0,1]\), определяемого тождеством
\[
	\forall y,z\in\Wo_2^1[0,1]\qquad\langle T_{\rho}(\lambda) y,z\rangle=
	\int\limits_0^1 \left\{y'\overline{z'}+\lambda P\cdot(y'\overline{z}+
	y\overline{z'})\right\}\,d\mu,
\]
где \(P\) "--- обобщённая первообразная весовой функции \(\rho\in
\Wo_2^{-1}[0,1]\). Как и в~\cite{VlSh}, именно такой пучок мы связываем с
граничной задачей~\eqref{eq:0.1},~\eqref{eq:0.2}. При формулировке дальнейших
утверждений будут использоваться обозначения и терминология из~\cite{VlSh}.

На основе теоремы~\ref{tm:1.1} устанавливается справедливость следующего
утверждения.
\begin{tm}\label{tm0}
Пусть \(P\in L_2[0,1]\) "--- арифметически самоподобная с шагом \(\nu\) функция,
имеющая положительный спектральный порядок \(D\). Пусть при этом справедливы
следующие условия:
\begin{enumerate}
\item обобщённая производная \(\rho\in\Wo_2^{-1}[0,1]\) функции \(P\) отлична
от нуля;
\item для любого номера \(k\leqslant n\) со свойством \(d_k>0\) отношение
\begin{equation}\label{eq:4.1}
	\dfrac{\ln\left(a_k\,|d_k|\right)}{\nu}
\end{equation}
является чётным;
\item для любого номера \(k\leqslant n\) со свойством \(d_k<0\)
отношение~\eqref{eq:4.1} является нечётным.
\end{enumerate}
Тогда существует такая непрерывная положительная
\mbox{\(2\)-пе}\-ри\-оди\-чес\-кая функция \(s\), что при \(\lambda\to+\infty\)
справедливо асимптотическое представление
\[
	\ind T_{\rho}(\lambda)=\lambda^{D/2}\cdot\left(s\left(\dfrac{\ln\lambda}%
	{\nu}\right)+o(1)\right),
\]
а при \(\lambda\to-\infty\) справедливо асимптотическое представление
\[
	\ind T_{\rho}(\lambda)=|\lambda|^{D/2}\cdot\left(s\left(
	\dfrac{\ln|\lambda|}{\nu}-1\right)+o(1)\right).
\]
\end{tm}

Аналогично, на основе теорем~\ref{tmA:1} и~\ref{tmA:2} устанавливается
справедливость следующего утверждения.
\begin{tm}\label{tm1}
Пусть \(P\in L_2[0,1]\) "--- неарифметически самоподобная функция, имеющая
положительный спектральный порядок \(D\). Пусть также \(\rho\in
\Wo_2^{-1}[0,1]\) "--- обобщённая производная функции \(P\). Тогда имеют место
следующие факты:
\begin{enumerate}
\item существуют такие неотрицательные числа \(s_{\pm}\), что при
\(\lambda\to\pm\infty\) справедливы асимптотические представления
\[
	\ind T_{\rho}(\lambda)=|\lambda|^{D/2}\cdot\left(s_{\pm}+o(1)\right);
\]
\item\label{ok:2} если при некотором \(k\in\{1,\ldots,n\}\) имеет место
неравенство \(d_k<0\), то справедливо равенство \(s_+=s_-\);
\item\label{ok:3} если для некоторой функции \(y\in\Wo_2^{-1}[0,1]\) имеет место
неравенство \(\langle\rho,|y|^2\rangle>0\), то число \(s_+\) положительно.
Аналогично, если для некоторой функции \(y\in\Wo_2^{-1}[0,1]\) имеет место
неравенство \(\langle\rho,|y|^2\rangle<0\), то число \(s_-\) положительно.
\end{enumerate}
\end{tm}

Доказательства теорем~\ref{tm0} и~\ref{tm1} полностью аналогичны доказательству
теоремы~4.2 из работы~\cite{VlSh}. Поэтому на деталях мы здесь не
останавливаемся.

\section{Численные результаты}
\begin{table}[t]
\begin{tabular}{|r|rrr|rrr|c|r|rrr|rrr|}
\cline{1-7} \cline{9-15}
{\(n\)}&\multicolumn{3}{|c|}{\(\lambda_n\)}&
\multicolumn{3}{|c|}{\(n/\lambda_n^{\log_6 2}\)}&
\hphantom{!}\hspace{1em}\hphantom{!}&{\(n\)}&
\multicolumn{3}{|c|}{\(\lambda_n\)}&
\multicolumn{3}{|c|}{\(|n|/|\lambda_n|^{\log_6 2}\)}\\ \cline{1-7} \cline{9-15}
1&\(6,72\cdot 10^0\)&\(\pm\)&\(1\%\)&0,478&\(\pm\)&0,001&
&\(-1\)&\(-2,30\cdot 10^1\)&\(\pm\)&\(1\%\)&0,297&\(\pm\)&0,001\\
2&\(9,75\cdot 10^1\)&\(\pm\)&\(1\%\)&0,340&\(\pm\)&0,001&
&\(-2\)&\(-2,58\cdot 10^1\)&\(\pm\)&\(1\%\)&0,569&\(\pm\)&0,001\\
3&\(1,40\cdot 10^2\)&\(\pm\)&\(1\%\)&0,443&\(\pm\)&0,001&
&\(-3\)&\(-5,11\cdot 10^2\)&\(\pm\)&\(1\%\)&0,269&\(\pm\)&0,001\\
4&\(1,50\cdot 10^2\)&\(\pm\)&\(1\%\)&0,575&\(\pm\)&0,001&
&\(-4\)&\(-5,82\cdot 10^2\)&\(\pm\)&\(1\%\)&0,341&\(\pm\)&0,001\\
5&\(2,34\cdot 10^2\)&\(\pm\)&\(1\%\)&0,605&\(\pm\)&0,001&
&\(-5\)&\(-5,86\cdot 10^2\)&\(\pm\)&\(1\%\)&0,425&\(\pm\)&0,001\\
6&\(2,89\cdot 10^3\)&\(\pm\)&\(1\%\)&0,275&\(\pm\)&0,001&
&\(-6\)&\(-8,12\cdot 10^2\)&\(\pm\)&\(1\%\)&0,449&\(\pm\)&0,001\\
7&\(3,06\cdot 10^3\)&\(\pm\)&\(1\%\)&0,313&\(\pm\)&0,001&
&\(-7\)&\(-8,41\cdot 10^2\)&\(\pm\)&\(1\%\)&0,517&\(\pm\)&0,001\\
8&\(3,06\cdot 10^3\)&\(\pm\)&\(1\%\)&0,358&\(\pm\)&0,001&
&\(-8\)&\(-9,03\cdot 10^2\)&\(\pm\)&\(1\%\)&0,575&\(\pm\)&0,001\\
9&\(3,49\cdot 10^3\)&\(\pm\)&\(1\%\)&0,383&\(\pm\)&0,001&
&\(-9\)&\(-9,10\cdot 10^2\)&\(\pm\)&\(1\%\)&0,645&\(\pm\)&0,001\\
10&\(3,49\cdot 10^3\)&\(\pm\)&\(1\%\)&0,426&\(\pm\)&0,001&
&\(-10\)&\(-1,41\cdot 10^3\)&\(\pm\)&\(1\%\)&0,605&\(\pm\)&0,001\\
11&\(3,51\cdot 10^3\)&\(\pm\)&\(1\%\)&0,467&\(\pm\)&0,001&
&\(-11\)&\(-1,69\cdot 10^4\)&\(\pm\)&\(1\%\)&0,254&\(\pm\)&0,001\\
12&\(3,53\cdot 10^3\)&\(\pm\)&\(1\%\)&0,509&\(\pm\)&0,001&
&\(-12\)&\(-1,73\cdot 10^4\)&\(\pm\)&\(1\%\)&0,275&\(\pm\)&0,001\\
\cline{1-7} \cline{9-15}
\end{tabular}

\vspace{0.5cm}
\caption{Оценки первых собственных значений для случая \(N=3\),
\(a_1=a_2=a_3=1/3\), \(d_1=d_3=-1/2\), \(d_2=0\), \(\beta_1=0\),
\(\beta_2=\beta_3=1/2\).}
\label{tab:1}
\end{table}
В таблице~\ref{tab:1} представлены результаты численных расчётов для первых
двенадцати положительных и отрицательных собственных значений задачи
Штурма--Лиувилля, весовой функцией в которой выступает обобщённая производная
квадратично суммируемой функции с параметрами самоподобия
\begin{gather*}
	N=3,\\
	a_1=a_2=a_3=1/3,\\
	d_1=d_3=-1/2,\quad d_2=0,\\
	\beta_1=0,\quad \beta_2=\beta_3=1/2.
\end{gather*}
Эта задача подпадает под действие теоремы~\ref{tm0}. Отвечающая ей функция
\(s\) имеет в точках \(0\) и \(1\) значения
\begin{align*}
	0,29\leqslant &s(0)\leqslant 0,36,\\
	0,60\leqslant &s(1)\leqslant 0,68.
\end{align*}
Таким образом, утверждаемое теоремой~\ref{tm0} удвоение периода этой функции по
сравнению с рассмотренным в работе~\cite{VlSh} невырожденным случаем
представляет собой действительно наблюдаемый эффект.

\end{document}